\documentclass[12pt]{article}
\title{Reverse mathematics of the relativised fast growing hierarchy}
\author{Florian Pelupessy}
\date{}
\usepackage[english]{babel}
\usepackage{amsmath, amssymb}
\usepackage{url}
\usepackage{xcolor}
\usepackage{hyperref}
\usepackage{tikz}
\usepackage{multicol}
\usepackage{parskip}
\usepackage{makeidx}
  \makeindex

 \usepackage[initials]{amsrefs}

\BibSpec{article}{%
  +{}{\PrintAuthors}				{author} %
  +{,}{ \textit}					{title}
  +{,}{ }							{journal}
  +{}{ \textbf}					{volume}
  +{,}{ }							{note}
  +{}{ \parenthesize}				{date}
  +{:}{ }							{pages}
  +{.}{}							{transition}
  +{}{ }							{review}
  +{\\}{ \url }							{url}
}

\BibSpec{thesis}{%
  +{}{\PrintAuthors}				{author} %
  +{,}{ \textit}					{title}
  +{,}{ }							{journal}
  +{}{ \textbf}					{volume}
  +{,}{ }							{note}
  +{,}{ }							{date}
  +{:}{ }							{pages}
  +{.}{}							{transition}
  +{}{ }							{review}
  +{\\}{ \url }							{url}
}

\usepackage[T1]{fontenc}
\usepackage{libertine}

\newtheorem{theorem}{Theorem}
\newtheorem{lemma}{Lemma}
\newtheorem{definition}{Definition}

\newcommand{\qed}{
\begin{flushright}
$\Box$
\end{flushright}
}

\newcommand{\RCA}{\mathrm{RCA}_0}

\newcommand{\N}{\mathbb{N}}

\begin{document}
\maketitle
We examine the reverse mathematical status of the totality of the relativised fast growing hierarchy. We assume a suitable encoding for ordinals up to $\varepsilon_0$ as in e.g. \cite{hajekpudlak}, or alternatively as ordinal terms.  When we have transfinite recursion up to $\varepsilon_0$ available, we define a fast growing hierarchy relative to $f$ as follows:

\begin{tabular}{lcl}
$F^f_0 (x)$			& $=$ & $f(x)$, \\
$F^f_{\alpha+1} (x)$ 	& $=$ & ${F^f_{\alpha}}^{(x+1)}(1)$, \\
$F^f_\gamma (x)$	& $=$ & $F^f_{\gamma [x]} (x)$ \textrm{if $\gamma$ is a limit}. \\
\end{tabular}

The status of the totality of $F_\omega$, the relativised Ackermann function, has been determined in \cite{kreuzeryokoyama}. Unlike in that paper, we determine the status directly and for all $\alpha \leq \varepsilon_0$. Take $\omega_0=0$ and $\omega_{n+1}=\omega^{\omega_n}$. The main result of this note is:
\begin{theorem}
$\RCA$ proves that the following are equivalent for every $n$:
\begin{enumerate}
  \item $\omega_{n+1}$ is well founded: every strictly descending of ordinals below $\omega_{n+1}$ is finite. 
  \item $F^f_{\omega_n}$ is total for every $f\colon\N\rightarrow\N$.
\end{enumerate}
\end{theorem}

\subsection*{Definition}
Define the following function $K\colon (\varepsilon_0)^* \times \N \rightarrow (\varepsilon_0)^* \times \N$. Intuitively, this function represents one step in the obvious way to attempt to compute a value for $F^f$:
\[
K_f( \alpha_0 \dots \alpha_n,  x  )= \left\{ 
\begin{array}{ll}
(\alpha_0 \dots \alpha_{n-1}, f(x)) 													& \textrm{if $\alpha_n=0$}, \\
(\alpha_0 \dots \alpha_{n-1} \overbrace{\beta \dots \beta}^{\textrm{$x+1$ times}}, 1)	& \textrm{if $\alpha_n=\beta+1$}, \\
(\alpha_0 \dots \alpha_{n-1} \alpha_n [x], x)											& \textrm{if $\alpha_n$ is a limit}, \\
\end{array}
\right.
\]
and $K_f( \langle \rangle , x) =(\langle \rangle, x)$. Notice that $(\alpha_0 \dots \alpha_n, x)$ simply represents the term $F_{\alpha_0} ( \dots ( F_{\alpha_n} (x) \dots )$. $F^f$ is the result of repeated applications of the `computation steps' (when it exists).
\begin{definition} 
$F^f_\alpha (x) =\mu y. \exists nK_f^{(n)} (\alpha , x)= ( \langle \rangle , y )$. We call the sequence $\{K_f^{(i)} (\alpha , x)\}_{i\in \N}$ the derivation of $F^f_\alpha(x)$. 
\end{definition}
One can show that this definition is equivalent to usual $\Delta^0_1$ definitions as in, e.g. \cite{hajekpudlak} (adapted to take into account the different initial function and slightly different conditions). 
\subsection*{Strength}
\begin{lemma}[$\RCA$]
For $\alpha\leq\varepsilon_0$: $\omega^\alpha$ is well founded if and only if $F^f_\alpha$ is a total function for every $f\colon\N\rightarrow\N$. 
\end{lemma}
\emph{Proof:} ``$\Rightarrow$'': Take $h(\alpha_0 \dots \alpha_n, x)=\omega^{\alpha_0} + \dots + \omega^{\alpha_n}$ and $h(\langle \rangle , x)=0$. By well-foundedness the sequence $\{ h( K_f^{(i)} (\alpha , x))\}_{i \in\N}$ reaches zero.

\begin{definition}[Maximal coefficient]
$\mathrm{mc}(0)=0$ and, for $\alpha=\omega^{\alpha_0}\cdot a_0 + \dots + \omega^{\alpha_n}\cdot a_n$ with $\alpha_0 > \dots > \alpha_n$, $a_i >0$:
\[
\mathrm{mc}(\alpha)=\max\{ \mathrm{mc}(\alpha_i), a_i \}.
\]
\end{definition}

``$\Leftarrow$'': Given infinite sequence $\omega^\alpha=\alpha_0>\alpha_1> \alpha_2> \dots$, take $f(x)>\mathrm{mc}(\alpha_{x+1})+x+1$ and strictly increasing. We show that this implies that for every $i>0$ we have $h(K^{(i)}_f (\alpha, f(0))) > \alpha_i$, in contradiction with the totality of $F^f_\alpha$. 

First, notice:
\begin{enumerate}
  \item By $\Pi^0_1$-induction on $d$: if $\omega_d > \gamma > \beta$ and $\gamma$ is a limit, then $\gamma[\mathrm{mc}(\beta)+1] > \beta$.
  \item $F^f_\beta (y) \geq f(y)$, hence ${F^f}_\beta^{(y)} (1)> y$ for all $\beta\leq\alpha$, $y$ which occur in the derivation of $F^f_\alpha(f(0))$.
  \item  if $n$ is the smallest such that $K_f^{(n)} (\beta , y)= ( \langle \rangle , z )$, then $K_f^{(n)} (\sigma\beta , y)= ( \sigma , z )$.
\end{enumerate}
Notation: $K^{(i)}_j=(K^{(i)}_f (\alpha, f(0)))_j$ and $K^{(i)}=K^{(i)}_f (\alpha, f(0))$.

Take $a_0=0$ and $a_{i+1}=a_i+1$ if  $K^{(a_i)}_0$ ends with a zero, otherwise as follows:

Let $b\geq a_i$ be the smallest such that $K^{(b)}_0$ ends with a successor $\beta+1$, take:
\[
a_{i+1} = \min b+n+1 \textrm{ such that } K_f^{(n)} (\overbrace{\beta \dots \beta}^{i+1} , 1)=(\langle\rangle, z).
\]
\emph{Claim:} For every $i$ we have:
\[
h(K^{(a_i)}) > \alpha_{i+1}
\]
and
\[
K^{(a_i)}_1 \geq f(i).
\]
\emph{Proof of the claim:} Induction on $i$, if $i=0$ the claim follows directly. For the induction step, assume that the claim is true for $a_i$.

Case 1) $a_{i+1}=a_i+1$: The inequalities follow directly from the definition of $K$:
\[
h(K^{(a_{i+1})})=h(K^{(a_{i})})-1\geq \alpha_{i+1}>\alpha_{i+2}
\]
and 
\[
K^{(a_{i+1})}_1= f(K^{(a_{i})}_1) \geq f(f(i)) \geq f(i+1).
\]

Case 2) Let $b$ and $\beta$ be those from the definition of $a_{i+1}$. $K^{(j)}_0$ ends with a limit for $j \in [a_i, b)$ (if $b>a_i)$, hence, by induction hypothesis and notice (1), $h(K^{(b)}) > \alpha_{i+1}$ and $K^{(b)}_1 \geq f(i)$.  

$K^{(b)}_0$ is of the form $\gamma_0 \dots \gamma_l \beta+1$. Therefore, $K^{(b+1)}_0$ has the form:
\[
\gamma_1 \dots \gamma_l \overbrace{\beta \dots \beta}^{\geq \mathrm{mc}(\alpha_{i+1})+1} \overbrace{\beta \dots \beta}^{i+1},
\]
so $h(K^{(a_{i+1})})\geq h(K^{(b)})[\mathrm{mc}(\alpha_{i+1})+1] > \alpha_{i+1}> \alpha_{i+2}$ by notice (1) and (3).

By notice (2) and (3),  $K^{(a_{i+1})}_1 \geq {F^f}_\beta^{(i+1)}(1) \geq f(i+1)$. 

This ends the proof of the claim, hence the lemma.
\qed

\end{document}